\documentclass[12pt]{article}
\usepackage{url,color}
\usepackage{hyperref}
\hypersetup{urlcolor=red,colorlinks=true,pdfwindowui=true,pdffitwindow=true,pdfview=XYZ}
\hypersetup{bookmarksopen=true,bookmarksnumbered=true}
\hypersetup{pdftoolbar=true,pdfmenubar=true}
\begin{document}
\title{\huge\bf Toy operads
   \footnote{Supported by \href{http://www.etf.ee/index.php?keel=ENG}{EstSF} grant 5634}
   }
\date{}
\author{
\href{http://www.staff.ttu.ee/~eugen}{\Large Eugen Paal}\\
\href{http://www.ttu.ee/?lang=en}{Tallinn University of Technology}
\\
Ehitajate tee 5,
19086 \href{http://www.tallinn.ee/eng}{Tallinn},
\href{http://www.riik.ee/en}{Estonia}\\
\href{http://www.staff.ttu.ee/~eugen/gen/message.html}{eugen@edu.ttu.ee}
}
\maketitle
\thispagestyle{empty}
\vskip.5cm
\begin{center}
{\bf Abstract}
\end{center}
Didactic operadic entertainment for pedestrians based on prosper overlays, 4 slides.
The following visual toys are included.
\begin{itemize}
\itemsep-3pt
\item
Operad of little squares.
\item
Operad of planar rooted trees.
\item
Operad algebra example.
\end{itemize}
\par
\noindent
\href{http://www.staff.ttu.ee/~eugen/toy/toy_operad.pdf}{Read or download.}
Preferable browser is
\href{http://www.adobe.com/products/acrobat/readstep2.html}{Adobe Reader 7}.
\vskip1cm
\noindent
{\bf Key words.} Operad, internet, clicking, tree, grafting.
\par
\vskip0.5cm
\noindent
\href{http://www.ams.org}{\bf AMS}
\href{http://www.ams.org/msc}{\bf MSC}.
\href{http://www.ams.org/msc/18Dxx.html}{18D50}
\end{document}